\documentclass[12pt, reqno]{amsart}
\usepackage{amsmath, amsthm, amscd, amsfonts, amssymb, graphicx, color}
\usepackage[english]{babel}
\usepackage{cite}
\usepackage[bookmarksnumbered, colorlinks, plainpages]{hyperref}

\newtheorem{theorem}{Theorem}[section]
\newtheorem{lemma}[theorem]{Lemma}

\newtheorem{corollary}[theorem]{Corollary}
\theoremstyle{definition}

\newtheorem{problem}[theorem]{Problem}
\theoremstyle{remark}

\numberwithin{equation}{section}
\begin{document}
\setcounter{page}{1}

\title{Local and $2$-Local derivations on noncommutative Arens algebras}

\author[Sh.~A.~Ayupov, K.~K.~Kudaybergenov, B.~O.~Nurjanov, A.~K.~Alauatdinov]{Sh.~A.~Ayupov$^1,$
 K.~K.~Kudaybergenov$^2,$ B.~O.~Nurjanov$^3,$  and A.~K.~Alauatdinov$^4$}

\address{$^{1}$
Institute of
 Mathematics and Information  Technologies
 \\ Uzbekistan Academy of Sciences  \\
 100125  Tashkent,   Uzbekistan\\
 and
 the Abdus Salam International Centre\\
 for Theoretical Physics (ICTP),
  Trieste, Italy
\newline}

\email{\textcolor[rgb]{0.00,0.00,0.84}{sh$_{-}$ayupov@mail.ru}}

\address{$^{2}$ Department of Mathematics, Karakalpak state university\\
Ch. Abdirov 1,  230113, Nukus,    Uzbekistan
\newline}

\email{\textcolor[rgb]{0.00,0.00,0.84}{karim2006@mail.ru}}

\address{$^{3}$ Department of Mathematics, Karakalpak state university\\
Ch. Abdirov 1,  230113, Nukus,    Uzbekistan
\newline}

\email{\textcolor[rgb]{0.00,0.00,0.84}{nurjanov@list.ru}}

\address{$^{4}$ Institute of
 Mathematics and Information  Technologies
 \\ Uzbekistan Academy of Sciences  \\
 100125  Tashkent,   Uzbekistan\\}

\email{\textcolor[rgb]{0.00,0.00,0.84}{amir$_{-}$t85@mail.ru}}






\begin{abstract}
The paper is devoted to so-called local and $2$-local derivations on the noncommutative Arens  algebra
$L^\omega(M, \tau)$
associated with a  von Neumann algebra $M$ and a faithful normal
semi-finite trace $\tau.$ We prove that every $2$-local derivation on
$L^\omega(M, \tau)$ is a spatial derivation, and if $M$ is a finite von Neumann algebra, then each local
derivation on $L^\omega(M, \tau)$ is also a spatial derivation and every $2$-local derivation on $M$ is in fact an inner derivation.
\end{abstract}

\maketitle

\section{Introduction}

Given an algebra $\mathcal{A},$ a linear operator $D:\mathcal{A}\rightarrow \mathcal{A}$ is
called a \textit{derivation}, if $D(xy)=D(x)y+xD(y)$ for all $x,
y\in \mathcal{A}$ (the Leibniz rule). Each element $a\in \mathcal{A}$ implements a
derivation $D_a$ on $\mathcal{A}$ defined as $D_a(x)=[a, x]=ax-xa,\,x\in \mathcal{A}.$ Such
derivations $D_a$ are said to be \textit{inner derivations}. If
the element $a,$ implementing the derivation $D_a,$ belongs to a
larger algebra $\mathcal{B}$ containing $\mathcal{A},$ then $D_a$ is called \textit{a
spatial derivation} on $\mathcal{A}.$

In the papers \cite{Alb3}, \cite{Ayup}  S.~Albeverio and the first two authors have proved the spatiality of
derivations on the non commutative Arens algebra
$L^{\omega}(M, \tau)$ associated with an arbitrary von
Neumann algebra $M$ and a faithful normal semi-finite
trace $\tau.$ Moreover if the trace $\tau$ is finite then every
derivation on $L^{\omega}(M, \tau)$ is inner.

There exist various types of linear operators which are close to derivations
\cite{Bre3, Bre1, John, Kad, Lar, Kim, Lin},\cite{Sem1}. In particular R.~Kadison \cite{Kad} has introduced and
investigated so-called local derivations on von Neumann algebras and some polynomial algebras.

A linear operator $\Delta$ on an algebra $\mathcal{A}$ is called a
\textit{local derivation} if given any $x\in \mathcal{A}$ there exists a
derivation $D$ (depending on $x$) such that $\Delta(x)=D(x).$  The
main problems concerning this notion are to find conditions under
which local derivations become derivations and to present examples of algebras
with local derivations that are not derivations \cite{Kad}, \cite{Lar}. In particular
Kadison \cite{Kad} has proved that each
continuous local derivation from a von Neumann algebra
$M$ into a dual $M$-bimodule is a derivation.

Later this result has been extended in \cite{Bre3} to a larger class of
linear operators $\Delta$ from $M$ into a normed
$M$-bimodule $E$ satisfying the identity
$\Delta(p)=\Delta(p)p + p\Delta(p)$
for every idempotent $p\in M.$
In \cite{Bre1} it was proved that every linear operator $\Delta$
on the algebra $M_{n}(R)$ satisfying above identity  is automatically a
derivation, where $M_{n}(R)$ is the algebra of $n\times n$
matrices over a unital commutative ring $R$ containing $1/2.$

In \cite{John}, B.~E.~Johnson has extended Kadison's result and
proved that every local derivation from a $C^{\ast}$-algebra $\mathcal{A}$ into any Banach $\mathcal{A}$-bimodule
is a derivation. He also showed that every local derivation from
 a $C^{\ast}$-algebra $\mathcal{A}$ into any Banach $\mathcal{A}$-bimodule is continuous.
In \cite{Nur}, \cite{Nur1} local derivations have been investigated on the algebra $S(M)$
of all measurable operators  with respect a von Neumann algebra $M$. In particular, we have  proved
 that for type I  von Neumann algebras without abelian direct summands every local derivation
 on $S(M)$ is a derivation. Moreover, in the  case of abelian von Neumann algebra $M$ we obtained necessary
and sufficient conditions for the algebra $S(M)$ to admit local derivations that
are not derivations.

In 1997, P. Semrl \cite{Sem1}  introduced the concept of
$2$-local derivations.
A  map $\Delta:\mathcal{A}\rightarrow\mathcal{A}$  (not linear in general) is called a
 $2$-\emph{local derivation} if  for every $x, y\in \mathcal{A},$  there exists
 a derivation $D_{x, y}:\mathcal{A}\rightarrow\mathcal{A}$
such that $\Delta(x)=D_{x, y}(x)$  and $\Delta(y)=D_{x, y}(y).$
Local and $2$-local maps have been studied on different operator algebras by many
authors
\cite{Bre3,  Bre1, John,  Kad, Kim, Lar, Lin, Mol, Sem1}.

In \cite{Sem1}, P. Semrl described
$2$-local derivations on the algebra $B(H)$ of
 all bounded linear operators on the infinite-dimensional
separable Hilbert space $H.$  A similar
description for the finite-dimensional case appeared later in \cite{Kim},
\cite{Mol}. In the paper \cite{Lin}  $2$-local derivations have been described on matrix algebras over finite-dimensional division rings.

The present paper  is devoted to study of local and $2$-local derivations on the noncommutative Arens  algebra
$L^\omega(M, \tau)$
associated with a von Neumann algebra $M$ and a faithful normal
semi-finite trace $\tau.$

In section $2$ we show that if $M$ is a finite von Neumann algebra, then each local
derivation on $L^{\omega}(M, \tau)$ is in fact a spatial derivation. In section $3$  we consider
  $L^{\omega}(M, \tau)$ for an arbitrary von Neumann algebra  $M$  with a faithful normal
 semi-finite trace $\tau$ and prove that every $2$-local derivation on $L^{\omega}(M, \tau)$
 is also a spatial derivation. We also extend the mentioned results of \cite{Kim},\cite{Mol} and \cite{Sem1} to arbitrary finite von Neumann algebras.

\section{Local derivation on  Arens algebras}

Let  $B(H)$ be the $\ast$-algebra of all
bounded linear operators on a Hilbert space $H,$ and let $\textbf{1}$ be the identity operator
on $H.$ Consider a von Neumann algebra $M\subset B(H)$
 with a faithful normal semi-finite trace
$\tau.$ Denote by $P(M)=\{p\in M: p=p^2=p^\ast\}$ the lattice of all projections in
$M.$

A linear subspace  $\mathcal{D}$ in  $H$ is said to be
\emph{affiliated} with  $M$ (denoted as  $\mathcal{D}\eta M$), if
$u(\mathcal{D})\subset \mathcal{D}$ for every unitary  $u$ from
the commutant
$$M'=\{y\in B(H):xy=yx, \,\forall x\in M\}$$ of the von Neumann algebra $M.$

A linear operator  $x: \mathcal{D}(x)\rightarrow H,$ where the  domain  $\mathcal{D}(x)$
of $x$ is a linear subspace of $H,$ is said to be \textit{affiliated} with  $M$ (denoted as  $x\eta M$) if
$\mathcal{D}(x)\eta M$ and $u(x(\xi))=x(u(\xi))$
 for all  $\xi\in
\mathcal{D}(x)$  and for every unitary  $u\in M'.$

A linear subspace $\mathcal{D}$ in $H$ is said to be \textit{strongly
dense} in  $H$ with respect to the von Neumann algebra  $M,$ if

1) $\mathcal{D}\eta M;$

2) there exists a sequence of projections
$\{p_n\}_{n=1}^{\infty}$ in $P(M)$  such that
$p_n\uparrow\textbf{1},$ $p_n(H)\subset \mathcal{D}$ and
$p^{\perp}_n=\textbf{1}-p_n$ is finite in  $M$ for all
$n\in\mathbb{N}$.

A closed linear operator  $x$ acting in the Hilbert space $H$ is said to be
\textit{measurable} with respect to the von Neumann algebra  $M,$ if
 $x\eta M$ and $\mathcal{D}(x)$ is strongly dense in  $H.$

 Denote by $S(M)$  the set of all linear operators on $H,$ measurable with
respect to the von Neumann algebra $M.$ If $x\in S(M),$
$\lambda\in\mathbb{C},$ where $\mathbb{C}$  is the field of
complex numbers, then $\lambda x\in S(M)$  and the operator
$x^\ast,$  adjoint to $x,$  is also measurable with respect to $M$
(see \cite{Seg}). Moreover, if $x, y \in S(M),$ then the operators
$x+y$  and $xy$  are defined on dense subspaces and admit closures
that are called, correspondingly, the strong sum and the strong
product of the operators $x$  and $y,$  and are denoted by
$x\stackrel{.}+y$ and $x \ast y.$ It was shown in  \cite{Seg} that
$x\stackrel{.}+y$ and $x \ast y$ belong to $S(M)$ and
these algebraic operations make $S(M)$ a $\ast$-algebra with the
identity $\textbf{1}$  over the field $\mathbb{C}.$ Here, $M$ is a
$\ast$-subalgebra of $S(M).$ In what follows, the strong sum and
the strong product of operators $x$ and $y$  will be denoted in
the same way as the usual operations, by $x+y$  and $x y.$

A closed linear operator $x$ in  $H$  is said to be \emph{locally
measurable} with respect to the von Neumann algebra $M,$ if $x\eta
M$ and there exists a sequence $\{z_n\}_{n=1}^{\infty}$ of central
projections in $M$ such that $z_n\uparrow\textbf{1}$ and $z_nx \in
S(M)$ for all $n\in\mathbb{N}$ (see \cite{Yea}).

Denote by $LS(M)$ the set of all linear operators that are locally
measurable with respect to $M.$ It was proved in \cite{Yea}  that
$LS(M)$ is a $\ast$-algebra over the field $\mathbb{C}$ with
identity $\textbf{1},$ the operations of strong addition, strong
multiplication, and passing to the adjoint. In such a case, $S(M)$
is a $\ast$-subalgebra in $LS(M).$ In the case where $M$ is a
finite von Neumann algebra or a factor, the algebras $S(M)$ and
$LS(M)$ coincide. This is not true in the general case. In
\cite{Mur2} the class of von Neumann algebras $M$ has been described for
which the algebras  $LS(M)$ and  $S(M)$ coincide.

 Let   $\tau$ be a faithful normal semi-finite trace on $M.$ We recall that a closed linear operator
  $x$ is said to be  $\tau$\textit{-measurable} with respect to the von Neumann algebra
   $M,$ if  $x\eta M$ and   $\mathcal{D}(x)$ is
  $\tau$-dense in  $H,$ i.e. $\mathcal{D}(x)\eta M$ and given   $\varepsilon>0$
  there exists a projection   $p\in M$ such that   $p(H)\subset\mathcal{D}(x)$
  and $\tau(p^{\perp})<\varepsilon.$
   Denote by  $S(M,\tau)$ the set of all   $\tau$-measurable operators with respect to  $M.$

    It is well-known that $S(M)$ and $S(M, \tau)$ are $\ast$-subalgebras in $LS(M)$ (see \cite{Mur}).

Given $p\geq1$ put $L^{p}(M, \tau)=\{x\in S(M,
\tau):\tau(|x|^{p})<\infty\}.$ It is known   \cite{Mur} that  $L^{p}(M, \tau)$ is a Banach
space with respect to the norm
$$\|x\|_p=(\tau(|x|^{p}))^{1/p},\quad x\in L^{p}(M, \tau).$$
  Consider the intersection
\begin{center}
$L^{\omega}(M, \tau)=\bigcap\limits_{p\geq1}L^{p}(M, \tau).$
\end{center}
It is proved in \cite{Abd} that  $L^{\omega}(M, \tau)$ is a locally convex
complete metrizable  $\ast$-algebra with respect to the  topology
$t$ generated by the family of norms $\{\|\cdot\|_p\}_{p\geq1}.$
 The algebra
$L^{\omega}(M, \tau)$ is called a (non commutative) \textit{Arens
algebra}.

Note that $L^{\omega}(M, \tau)$ is  $\ast$-subalgebra in $S(M, \tau)$ and if $\tau$ is a finite trace then $M\subset L^\omega(M, \tau).$

Further  consider the following spaces
$$L^{\omega}_2(M, \tau)=\bigcap\limits_{p\geq 2}L^p(M, \tau)$$
and
$$
M+
L^{\omega}_2(M, \tau)=\{x+y: x\in M,
y\in L^{\omega}_2(M, \tau)\}.
$$

It is known \cite{Alb3} that
$L^{\omega}_2(M, \tau)$  and $M+ L^{\omega}_2(M,
\tau)$ are a
$\ast$-algebras and $L^{\omega}(M, \tau)$ is an ideal in $M+L^{\omega}_2(M,
\tau).$

Note that if $\tau(\textbf{1})<\infty$ then $M+ L^{\omega}_2(M,
\tau)=L^{\omega}_2(M,
\tau)=L^{\omega}(M, \tau).$

It is known \cite[Theorem 3.7]{Alb3} that if
   $M$ is a  von Neumann algebra with a faithful normal
semi-finite trace $\tau$ then any derivation $D$ on
  $L^{\omega}(M, \tau)$ is spatial, moreover it is implemented by an element of $M+L^{\omega}_2(M, \tau),$ i. e.
\begin{equation}
\label{SP}
D(x)=ax-xa, \quad x\in L^{\omega}(M, \tau),
\end{equation}
for some
$a\in M+L^{\omega}_2(M, \tau).$

Now we need several lemmata.

\begin{lemma}\label{A}
Let  $a\in LS(M), p_n\in P(M), n\in \mathbb{N},$  $p_n\uparrow \textbf{1}$
and  $p_n a p_n=0$ for all   $n\in \mathbb{N}.$ Then  $a=0.$
\end{lemma}

\begin{proof} Fix $k\in \mathbb{N}$ and take any  $n\geq k.$
Then  $p_k=p_k p_n$ and therefore from  $p_n a p_n=0$ it follows that
$p_k a p_n=p_k (p_n a p_n)=0.$ Thus  $(p_k a) p_n (p_k a)^\ast=0.$
Since  $p_n\uparrow \textbf{1},$ then  $(p_k a)(p_k a)^\ast=0,$ i.e.
$p_k a=0.$ Thus  $a^\ast p_k a=0.$
Again from   $p_k\uparrow \textbf{1},$ we obtain that  $a^\ast a=0,$ i.e.
$a=0.$ The proof is complete.
\end{proof}

Recall that the subalgebra $\mathcal{A}\subset LS(M)$ is said to be solid,
 if $x\in \mathcal{A},\, y\in LS(M),\,|y|\leq |x|$
implies $y\in \mathcal{A}.$

\begin{lemma}\label{B}
Let $M$ be a von Neumann algebra and  let $\mathcal{A}\subseteq LS(M)$ be a solid
$\ast$-subalgebra such that $M\subseteq \mathcal{A}.$ If   $\Delta: \mathcal{A}
\rightarrow \mathcal{A}$ is a local derivation such that
 $\Delta|_M\equiv 0,$ then $\Delta\equiv 0.$
 \end{lemma}

\begin{proof} Let    $\Delta: \mathcal{A}
\rightarrow \mathcal{A}$ be a local derivation such that
 $\Delta|_M\equiv 0.$ First let us  show that
\begin{equation}
\label{E}
 e\Delta(e^\perp x e^\perp) e=0
\end{equation}
 for every projection $e\in M$ and each element $x\in \mathcal{A}.$
 Indeed take a derivation $D: \mathcal{A}
\rightarrow \mathcal{A}$ such that
$$
 \Delta(e^\perp x e^\perp)=D(e^\perp x e^\perp).
$$
Then
$$
 e\Delta(e^\perp x e^\perp) e=
  e D(e^\perp x e^\perp) e=
  $$
  $$
  =e[D(e^\perp) x e^\perp + e^\perp D(x) e^\perp +e^\perp x D(e^\perp)]e=0,
$$
i.e.  $e\Delta(e^\perp x e^\perp) e=0.$

Now we take any positive element $x\in \mathcal{A}.$
 Put  $p_n=e_n(x),$ where  $\{e_\lambda(x)\}_{\lambda\in \mathbb{R}}$ is the spectral family
 of the element
 $x.$
 Since  $x p_n\in M,$
 then  $\Delta(xp_n)=0$ for all  $n\in \mathbb{N}.$  Therefore
   $$
\Delta(x)=\Delta(x p_n+x p_n^\perp)=
\Delta(x p_n)+\Delta(x p_n^\perp)=
\Delta(x p_n^\perp)=\Delta(p_n^\perp x p_n^\perp),
$$
i.e.
$$
\Delta(x)=\Delta(p_n^\perp x p_n^\perp).
$$
Thus by the equality (\ref{E}) we obtain that
 $$
p_n\Delta(x)p_n=p_n \Delta(p_n^\perp x p_n^\perp)  p_n=0.
$$
i.e.
 $$
p_n\Delta(x)p_n=0.
$$
Since   $p_n\uparrow \textbf{1},$  by  Lemma \ref{A} we obtain that
 $$
\Delta(x)=0
$$
for all  $x\geq 0.$ Since $\mathcal{A}$ is a solid $\ast$-subalgebra in $LS(M)$
 any element  $x$ from  $\mathcal{A}$ can be a represented
as a finite linear combinations of positive elements from  $\mathcal{A},$ therefore
$\Delta(x)=0$ for all $x\in \mathcal{A}.$
The proof is complete.
\end{proof}

\begin{lemma} \label{C}
Let $M$ be a von Neumann algebra with the center $Z(M)$ and  let $\mathcal{A}\subseteq LS(M)$ be a
$\ast$-subalgebra such that $M\subseteq \mathcal{A}.$
Then every local derivation $\Delta$ on the algebra
$\mathcal{A}$ is necessarily $P(Z(M))$-homogeneous, i.e.
$$\Delta(zx)=z\Delta(x)$$ for any central projection $z\in
P(Z(M))=P(\mathcal{M})\cap Z(M)$ and for all $x\in
\mathcal{A}.$
\end{lemma}

\begin{proof}
Take $z\in P(Z(M))$ and $x\in \mathcal{A}.$ For the element
$zx$ by the definition of the local derivation $\Delta$ there
exists a derivation $D$ on $\mathcal{A}$ such that
$\Delta(zx)=D(zx).$ Since the projection $z$ is central, one has
that  $D(z)=0$ and therefore
$$\Delta(zx)=D(zx)=D(z)x+zD(x)=zD(x),$$
 i.e. $\Delta(zx)=zD(x).$
  Multiplying this equality by $z$ we obtain
  $$z\Delta(zx)=z^2D(x)=zD(x)=\Delta(zx),$$  i.e.
$$z^{\bot}\Delta(zx)=(\mathbf{1}-z)\Delta(zx)=0.$$

Therefore by the linearity of $\Delta$ we have
$$
z\Delta(x)=z\Delta(zx)+z\Delta(z^{\bot}x)=z\Delta(zx)=\Delta(zx),
$$
and thus $z\Delta(x)=\Delta(zx).$ The proof is complete.
\end{proof}

The following theorem is the main result of this section.

\begin{theorem}\label{T1}
Let   $M$ be a  finite von Neumann algebra  with a faithful normal
semi-finite trace $\tau.$ Then any local derivation
$\Delta$ on the algebra  $L^{\omega}(M,\tau)$ is a spatial derivation of the form
(\ref{SP}).
\end{theorem}

\begin{proof} First assume that $\tau$ is a finite trace. Then
$M\subset L^\omega(M, \tau).$
Since  $L^{p}(M,\tau)$ ($p\geq 1$) is a  $M$-bimodule, then by Johnson's theorem
 (see  \cite[Theorem 7.5]{John}) the restriction of the local derivation  $\Delta$ on  $M$
is continuous as a map from  $M$ into  $L^{p}(M,\tau)$ for all  $p\geq 1.$ From
 \cite[Theorem 2.3]{Bre1} it follows that the
continuous local derivation  $\Delta|_M:M\rightarrow L^{p}(M,\tau)$
is a derivation.  Now by  \cite[Lemma]{Ayup}
the derivation  $\Delta|_M : M\rightarrow L^{\omega}(M,\tau)$ is spatial, i.e.
\begin{equation}\label{F}
\Delta(x)=ax-xa, \, x\in M
\end{equation}
for some $a\in L^{\omega}(M,\tau).$

Let us show that  $\Delta(x)=ax-xa$ for all  $x\in L^{\omega}(M,\tau).$
Consider the local derivation  $\Delta_0=\Delta-D_a.$
 Then from  (\ref{F}) we obtain that  $\Delta_0|_M\equiv 0.$ By Lemma~\ref{B} it follows that $\Delta_0\equiv 0.$
 This means that  $\Delta=D_a.$

Now let $\tau$ be a semi-finite trace. First note that if $M$ is a finite von Neumann algebra with a faithful normal semi-finite trace $\tau$, then the restriction of $\tau$ to the center $Z(M)$ of $M$ is also a semi-finite trace on $Z(M)$ (see for details the remark before Proposition 1.2 on page 2921 of \cite {Alb2}).  Therefore there exists a family of
orthogonal central projections  $\{z_i\}_{i\in I}$ in $M$ with
$\sup\limits_{i\in I}z_{i}=\textbf{1},$ such that $\tau(z_i)<\infty$
for all $i \in I.$ By Lemma \ref{C} we have that
$$\Delta(z_{i}x)=z_{i}\Delta(x)$$ for each $i\in I.$ This implies that $\Delta$ maps each
$z_{i}L^{\omega}(M,\tau)$ into
itself and hence induces a local derivation
$\Delta_{i}=\Delta|_{z_{i}L^{\omega}(M,\tau)}$ on the algebra
$z_i L^{\omega}(M,\tau)\cong
 L^{\omega}(z_i M,\tau_i)$ for each $i\in I,$ where $\tau_i$ is the restriction of $\tau$ on
$z_i M.$ Since $\tau(z_i)<\infty$  from the above it follows that the operator
$\Delta_{i}$  is a derivation.
Therefore for $x, y \in L^{\omega}(M,\tau)$ we have that
$$
\Delta_i((z_i x)(z_i y))=\Delta_{i}(z_i x) z_i y+z_i x\Delta_i(z_i y)
$$
for all $i\in I.$
By Lemma  \ref{C} we obtain that
$$
z_i \Delta(xy)=z_i[\Delta(x)y+x\Delta(y)],
$$
i.e.
$$
\Delta(x y)=\Delta(x) y+ x\Delta(y).
$$
This means that the operator  $\Delta$  is also a derivation
on $L^{\omega}(M,\tau)$
and  by \cite[Theorem 3.7]{Alb3} it can be  represented in the
form (\ref{SP}).
The proof is complete.
\end{proof}

Theorem \ref{T1} implies the following

\begin{corollary}\label{C1}
Let   $M$ be a   commutative von Neumann algebra  with a faithful normal
semi-finite trace $\tau.$ Then every local derivation
$\Delta$ on the algebra  $L^{\omega}(M,\tau)$ is identically zero.
\end{corollary}

 We conclude the section with the following open problem.
 \begin{problem}\label{P1}
 Is Theorem \ref{T1} valid for an arbitrary semi-finite von Neumann algebra $M$?
 \end{problem}

\section{$2$-local derivation on  Arens algebras}

Let now $M$ be an arbitrary semi-finite von Neumann algebra with a faithful
normal semi-finite trace $\tau.$ Consider the non commutative Arens algebra
$L^{\omega}(M,\tau).$ In order to prove the main result of the present section we need following simple lemma.

\begin{lemma}\label{F}
The algebra $L^{\omega}(M, \tau)$ is semi-prime, i.e. if  $a\in L^{\omega}(M, \tau)$ and
$a L^{\omega}(M, \tau) a=\{0\}$
then $a=0.$
\end{lemma}

\begin{proof}
Let $a\in L^{\omega}(M, \tau)$ and $a L^{\omega}(M, \tau) a=\{0\},$
i.e. $axa=0$ for all $x\in L^{\omega}(M, \tau).$ In particular for
$x=a^{*}$ we have $aa^{*}a=0$ and hence $a^{*}aa^{*}a=0,$ i.e.
$|a|^4=0.$ Therefore $a=0.$ The proof is complete.
\end{proof}

The following theorem is the main result of this section.

\begin{theorem}\label{TL}
Let   $M$ be a   von Neumann algebra  with a faithful normal
semi-finite trace $\tau.$ Then any 2-local derivation
$\Delta$ on the algebra  $L^{\omega}(M,\tau)$ is a derivation and has the form
(\ref{SP}).
\end{theorem}

\begin{proof} Let $\Delta:L^{\omega}(M,\tau)\rightarrow L^{\omega}(M,\tau)$ be
 a $2$-local derivation. Following the method of \cite[Theorem 1.2]{Lin} let us prove the additivity of $\Delta$.
For each $x, y\in L^{\omega}(M,\tau)$ there exists a derivation $D_{x, y}$ on $L^{\omega}(M,\tau)$ such that
$\Delta(x)=D_{x, y}(x)$  and $\Delta(y)=D_{x, y}(y).$
By \cite[Theorem 3.7]{Alb3} there exists element $a\in M+L^{\omega}_2(M,\tau)$
such that
$$
[a, xy]=D_{x, y}(xy)=D_{x, y}(x)y+xD_{x, y}(y)=\Delta(x)y+x\Delta(y),
$$
i.e.
$$
[a, xy]=\Delta(x)y+x\Delta(y).
$$
Now recall that  $L^{\omega}(M,\tau) \subset L^{1}(M,\tau)$ and
hence the trace $\tau$ accepts finite values on the algebra
$L^{\omega}(M,\tau).$  Since $L^{\omega}(M, \tau)$ is an ideal in
$M+L^{\omega}_2(M, \tau)$ \cite [Proposition 2.6]{Alb3} we have
$$
\tau(axy) = \tau((ax)y) = \tau(y(ax)) = \tau((ya)x) = \tau(x(ya))= \tau(xya).
$$
Thus we have
$$
0 = \tau(axy-xya)=\tau([a, xy])=\tau\left(\Delta(x)y+x\Delta(y)\right),
$$
i.e.
$\tau(\Delta(x)y)=-\tau(x \Delta(y)).$
For arbitrary $u, v, w\in L^{\omega}(M,\tau),$
 set $x=u+v,$ $y=w.$ Then from above we obtain
$$
\tau(\Delta(u+v)w)=-\tau((u+v)\Delta(w))=
$$
$$
=-\tau(u\Delta(w))- \tau(v\Delta(w))=\tau(\Delta(u)w)+\tau(\Delta(v)w)=\tau((\Delta(u)+\Delta(v))w),
$$
and so
$$
\tau((\Delta(u+v)-\Delta(u)-\Delta(v))w)=0
$$ for all $u, v, w\in L^{\omega}(M,\tau).$
Denote
$
b=\Delta(u+v)-\Delta(u)-\Delta(v)
$ and put
$
w=b^\ast.$ Then
$\tau(bb^\ast)=0.$ Since $\tau$ is faithful one has $bb^\ast=0,$ i.e. $b=0.$
Therefore
$$
\Delta(u+v)=\Delta(u)+\Delta(v),
$$
i.e. $\Delta$ is an additive operator.

Now let us  show that $\Delta$ is  homogeneous. Indeed, for
 each $x\in L^{\omega}(M,\tau), \lambda\in \mathbb{C}$ there exists a derivation $D_{x, \lambda x}$
 such that $\Delta(x)=D_{x, \lambda x}(x)$ and $\Delta(\lambda x)=D_{x, \lambda x}(\lambda x).$ Then
$$
\Delta(\lambda x)=D_{x, \lambda x}(\lambda x)=\lambda
D_{x, \lambda x}(x)=\lambda\Delta(x).
$$
Hence, $\Delta$ is homogenous and therefore is a linear operator.

Finally, for each $x\in L^{\omega}(M,\tau),$ there exists a derivation $D_{x, x^2}$
 such that $\Delta(x)=D_{x, x^2}(x)$ and $\Delta(x^2)=D_{x, x^2}(x^2).$ Then
$$
\Delta(x^2)=D_{x, x^2}(x^2)=D_{x, x^2}(x)x+xD_{x, x^2}(x)=\Delta(x)x+x\Delta(x)
$$
 for all $x\in L^{\omega}(M,\tau).$
Therefore, $\Delta$ is a linear Jordan derivation on
$L^{\omega}(M,\tau)$ in the sense of \cite{Bre2}. In
\cite[Theorem 1]{Bre2} it is proved that any Jordan derivation on
a semi-prime algebra is a derivation. Thus Lemma \ref{F} implies that
the linear operator $\Delta$ is a derivation on
$L^{\omega}(M,\tau).$ Now
 by \cite[Theorem 3.7]{Alb3} the derivation $\Delta$   can be  represented in
form (\ref{SP}).
The proof is complete.
\end{proof}

Theorem \ref{TL} implies the following

\begin{corollary}\label{C1}
Let   $M$ be a   commutative von Neumann algebra  with a faithful normal
semi-finite trace $\tau.$ Then any 2-local derivation
$\Delta$ on the algebra  $L^{\omega}(M,\tau)$ is identically zero.
\end{corollary}

In fact the method of the proof of Theorem \ref{TL}  can be easily
modified to the case of 2-local derivation on von Neumann algebras
with a separating family of normal finite traces  $\{\tau_{i}: i\in I\}.$
Since each finite von Neumann algebra has a separating family of
normal finite traces and all derivations on von Neumann algebras
are inner, we obtain the following analogue of \cite[Theorem
2]{Sem1} and \cite[Theorem 3]{Kim} for finite von Neumann
algebras:

\begin{theorem}\label{T3}
Every $2$-local derivation on a finite von Neumann algebra
is an inner derivation.
\end{theorem}

The latter theorem  gives rise to the following problem.
 \begin{problem}\label{P2}
 Is Theorem \ref{T3}  valid for an arbitrary  von Neumann algebra $M$?
 \end{problem}

\section*{Acknowledgments}

The first and the second named authors would like to acknowledge
the hospitality of the "Institut f\"{u}r Angewandte Mathematik",
Universit\"{a}t Bonn (Germany). This work is supported in part by
the DFG AL 214/36-1 project (Germany). The second author would also like to acknowledge
 the support of the German Academic Exchange Service -- DAAD .

\end{document}